\documentclass[12pt]{article}
\newlength{\mytopmargin}
\newlength{\myleftmargin}
\usepackage{amsmath,amsthm,amssymb,amsbsy,epsfig,graphicx,color,multicol,subfigure}
\usepackage{array,calc}
\usepackage[enableskew]{youngtab}
\usepackage{rotating}
\usepackage{young,epic}
\usepackage{a4wide,bm}
\usepackage{url}
\usepackage{hyperref}

\newtheorem{prop}{Proposition}

\newtheorem{cor}{Corollary}

\usepackage{amsmath,amsfonts,amssymb}

\usepackage{graphicx}

\begin{document}
%
%\begin{frontmatter}

\title{Lyapunov exponents for products of complex Gaussian random matrices}
\author{Peter J. Forrester }
\date{}
\maketitle
\noindent
\thanks{\small Department of Mathematics and Statistics, 
The University of Melbourne,
Victoria 3010, Australia email:  p.forrester@ms.unimelb.edu.au 
}

\begin{abstract}
The exact value of the Lyapunov exponents for the random matrix product $P_N = A_N A_{N-1} \cdots A_1$ with
each $A_i =  \Sigma^{1/2} G_i^{\rm c}$, where $\Sigma$ is a fixed $d \times d$ positive definite matrix and $G_i^{\rm c}$ a $d \times d$ complex Gaussian matrix with entries standard
complex normals, are calculated. Also obtained is an exact expression for the sum of the Lyapunov exponents in both the complex and real
cases, and the Lyapunov exponents for diffusing complex matrices.
\end{abstract}

\section{Introduction}

 Presently there is a great deal of interest, from the viewpoints of probability theory and applications to physics and communications
engineering, in the statistical properties of large random matrices (see e.g.~the recent texts \cite{AGZ09,PS11,Fo10,TV04}). Looking back to the
mid 1980's it would have been fair to say that the same applied to the topic of products of $d \times d$ random matrices. Thus it was in that era that the
foundational probabilistic works of Kesten, Furstenberg, Oseledec and others in the 1960's and 70's had matured to the extent that a book on the
subject was written \cite{BL85}; that a summer research conference was held with this topic dominating the subsequent proceedings \cite{CKN86}; and that
a number of applications to physics were in the full swing of investigation, culminating in the appearance of the research monograph \cite{CPV93}.

There have been some present day works that have aimed to combine the contemporary interest in the eigenvalues of large
random matrices with the topic of products of random matrices, by studying eigenvalue distributions of products of random matrices,
in the limit that the size of the matrices is large
\cite{BJW10,BJLNS11,PZ11,Ka08,OS10,Tu11}. There has also been a good deal of present day activity relating to the numerical
computation of Lyapunov exponents \cite{WT01,Ba07,Fi08,Ba09,Va10,Po10}. Regarding the latter, let
\begin{equation}\label{PA}
P_N = A_N A_{N-1} \cdots A_1,
\end{equation}
where each $A_i$ is a $d \times d$ independent, identically distributed random matrix such that the diagonal elements of $A^\dagger A$ have
finite second moments. According to the multiplicative ergodic theorem of Oseledec \cite{Os68,Ra79}, one has that the limiting matrix
\begin{equation}\label{PA1}
V_d := \lim_{N \to \infty} (P_N^\dagger P_N)^{1/(2N)}
\end{equation}
is well defined, with $d$ positive real eigenvalues $e^{\mu_1} \ge e^{\mu_2} \ge \cdots \ge e^{\mu_d}$. The $\{\mu_i\}$ are referred to as the
Lyapunov exponents. Already in 1973 Kingman had nominated methods to compute the largest Lyapunov exponent  \cite{Ki73} as an
outstanding problem in the field.
The work \cite{Po10} solves this problem in the case when each $A_i$ is chosen from a finite set, and has positive entries.

Implicit in the need for efficient computational methods is that it is generally not possible to compute the Lyapunov exponents analytically.
Some noteworthy exceptions occur in the case $d=2$; see e.g.~\cite{Ma93,MTW08,CTT10} and references therein. For general $d$, apart from
the case of diagonal matrices, it seems that the only exact computation of the Lyapunov exponents recorded in the literature is when the
$A_i$ are real Gaussian matrices with entries independent standard real normals. Then it is a result of Newman \cite{Ne86} that
\begin{equation}\label{N}
\mu_i = {1 \over 2} \Big ( \log 2 + \Psi \Big ( {d - i +1 \over 2} \Big ) \Big ) \qquad (i=1,\dots,d),
\end{equation}
where $\Psi(x)$ denotes the digamma function.

It is the aim of this paper to extend available exact results on the evaluation of Lyapunov exponents for certain $d \times d$ random matrices.
In particular, a closed form evaluation of $\{\mu_i\}$ is obtained for
\begin{equation}\label{AG}
A_i = \Sigma^{1/2} G_i^{\rm c},
\end{equation}
where $\Sigma$ is a fixed $d \times d$ positive definite matrix and $G_i^{\rm c}$ a $d \times d$ complex Gaussian matrix with entries standard
complex normals. Partial results are also obtained relating to the exact evaluation of $\{\mu_i \}$ for $A_i = \Sigma^{1/2} G_i^{\rm r}$, where
$G_i^{\rm r}$ denotes a $d \times d$ real Gaussian with independent standard normals as entries. The latter supplements the
$\Sigma = I_d$ result (\ref{N}), and the $d=2$ result
\begin{equation}\label{5.1}
\mu_1 = - {1 \over 2} \gamma + {1 \over 2} \log \Big ( {1 \over 2} {\rm Tr} \, \Sigma + \sqrt{ \det \Sigma} \Big ),
\end{equation}
where $\gamma$ denotes Euler's constant.

Our first result is the analogue of (\ref{N}) for  complex Gaussian matrices with entries independent standard complex normals. 

\begin{prop}\label{P1}
Consider the matrix product (\ref{PA}), with each $A_i$ given by (\ref{AG}) with $\Sigma = I_d$. We have
\begin{equation}\label{r1}
\mu_i = \Psi(d-i+1).
\end{equation}
\end{prop}

This result can be generalized by allowing for general $\Sigma$ in (\ref{AG}).

\begin{prop}\label{P2}
Consider the matrix product (\ref{PA}), with each $A_i$ given by (\ref{AG}) for general positive definite  $\Sigma$. 
Let the eigenvalues of $\Sigma^{-1}$ be denoted $\{y_j\}_{j=1,\dots,d}$. We have
\begin{equation}\label{r2}
\mu_k = - {1 \over 2 \prod_{1 \le i < j \le d} (y_j - y_i)}
\det \begin{bmatrix} [y_j^{i-1}]_{i=1,\dots,k-1 \atop j =1,\dots,d} \\
[(\log y_j) y_j^{k-1}]_{j=1,\dots,d} \\
[y_j^{i-1}]_{i=k+1,\dots,d \atop j =1,\dots,d} \end{bmatrix} + {1 \over 2} \Psi(k).
\end{equation}
\end{prop}

\begin{cor}\label{C1}
In the setting of Proposition \ref{P2}, one has the sum rule
\begin{equation}\label{15}
\mu_1 + \mu_2 + \cdots + \mu_d ={1 \over 2}  \sum_{m=1}^d\Big (-  \log y_m + \Psi(m) \Big ).
\end{equation}
\end{cor}

We remark that for $d$ large, the fact that $\Psi(d) \sim \log d$ tells us that the eigenvalues $y_m$ must have
the scaled form $y_m/d =  Y(m/d)$, with $Y(0) = 0 $ a PDF on $(0,1)$, for (\ref{15}) to have a well defined average as $d \to \infty$.

The above results will be proved in the next section. It will furthermore be showed that Corollary \ref{C1} can be proved independent of
Proposition \ref{P2}, and this will allow an analogue of (\ref{15}) in the case of real random matrices $A_i = \Sigma^{1/2} G_i^{\rm r}$
to be derived (see eq.~(\ref{2.33}) below). We give too, in eq.~(\ref{LK1}),
the evaluation of the so-called generalized maximum Lyapunov exponent
\cite{CPV93}
\begin{equation}\label{Lq}
L(q) = \lim_{N \to \infty} {1 \over N} \log \Big \langle ||P_N||^q  \Big \rangle 
\end{equation}
in the case of (\ref{AG}). 
In section 3  we calculate the Lyapunov exponents for
\begin{equation}\label{CA}
A_i := \lim_{m \to \infty} \prod_{j=1}^m e^{C_j^{(i)}/m^{1/2}}, \qquad C_j^{(i)} = H_1 + i H_2,
\end{equation}
where $H_1$ and $H_2$ are Hermitian matrices distributed with density function proportional to $\exp(- {\rm Tr} H_1^2/2 w_1)$ and
$\exp(- {\rm Tr} H_2^2/2 w_2)$ respectively. 
Section 4 discusses features of $\mu_1$ as given by (\ref{r2}).

\section{Proofs}
\setcounter{equation}{0}
\subsection{Background theory}
A fundamental characterization of the Lyapunov exponents defined below (\ref{PA1}) is that they satisfy \cite{Os68,Ra79}
\begin{equation}\label{7.1}
\mu_1 + \cdots + \mu_k = \sup \lim_{N \to \infty} {1 \over N}  \log {\rm Vol}_k \{y_1(N),\dots,y_k(N) \} \qquad (k=1,\dots,d).
\end{equation}
In (\ref{7.1}) $y_j(N) := P_N y_j(0)$, the supremum is over all sets of linearly independent vectors
$\{y_1(0),\dots,y_k(0) \}$ and Vol${}_k$ refers to the (generalized) volume of the parallelogram generated by the given set of $k$ vectors.
In regards to the latter, with 
$$
B_N := [ y_1(N) \: y_2(N) \: \cdots  \: y_k(N)],
$$
so that $B_N$ is the $d \times k$ matrix with its columns given by the $k$ vectors $\{y_j(N)\}$, we have
\begin{equation}\label{7.2}
 {\rm Vol}_k \{y_1(N),\dots,y_k(N) \} = \det (B_N^\dagger B_N)^{1/2} = \det (B_0^\dagger P_N^\dagger P_N B_0)^{1/2}
\end{equation}

Following \cite{CN84,Ne86} the basic fact that makes the computation of (\ref{7.1}) tractable for matrices (\ref{AG}) is that the distribution
of the random vector
$ G_i^{\rm c} \vec{x} / |\vec{x}|$, $\vec{x} \ne \vec{0}$, is independent of $\vec{x}$. Thus with $\{y_1(0),\dots,y_k(0) \}$ a set of $k$ linearly
independent unit vectors,, and $E_{d \times k}$ denoting the $d \times k$ matrix with 1's in the diagonal positions of the $k$ rows, and 0's
elsewhere, we have
\begin{equation}\label{8z}
B_0^\dagger P_N^\dagger P_N B_0 \mathop{=}\limits_{\rm d} 
\prod_{j=1}^N E_{d \times k}^T G_j^{\rm c \, \dagger} \Sigma G_j^{\rm c}  E_{d \times k} =
\prod_{j=1}^N G_{j,k}^{\rm c \, \dagger} \Sigma  G_{j,k}^{\rm c},
\end{equation}
where $G_{j,k}^{\rm c}$ denotes the $d \times k$ matrix formed by the first $k$ columns of $G_j^{\rm c}$.
Substituting in (\ref{7.2}), then substituting the result in (\ref{7.1}), we see firstly that there is no longer any dependence on
$\{y_1(0),\dots,y_k(0) \}$, so the sup operation in (\ref{7.1}) is redundant. We are then left with the expression
\begin{equation}\label{8}
\mu_1 + \cdots + \mu_k =  \lim_{N \to \infty} \sum_{j=1}^N {1 \over N} \log \det \Big ( G_{j,k}^{\rm c \, \dagger} \Sigma  G_{j,k}^{\rm c} \Big )^{1/2}.
\end{equation}
But
each $G_j^{\rm c}$ independently belongs to the set of complex rectangular $d \times k$
Gaussian matrices ${\mathcal N}_{d \times k}^{\rm c}(0,1)$
in which each entry is a standard complex normal.  The law of large numbers tells us that the limit in (\ref{8}) can be evaluated as an
average over this set,
\begin{equation}\label{8.1}
\mu_1 + \cdots + \mu_k =   \Big \langle \log \det  \Big ( G_{k}^{\rm c \, \dagger} \Sigma  G_{k}^{\rm c} \Big )^{1/2} \Big 
\rangle_{{\mathcal N}_{d \times k}^{\rm c}(0,1)}.
\end{equation}

\subsection{Proof of Proposition \ref{P1} and Corollary \ref{C1}}
Proposition \ref{P1} relates to the case $\Sigma = I_d$. Now the set of matrices ${\mathcal N}_{d \times k}^{\rm c}(0,1)$ have probability
density function proportional to $e^{-{\rm Tr}( G_j^{\rm c \, \dagger} G_j^{\rm c})}$. Thus the average in (\ref{8.1}) is a function of $ G_j^{\rm c \dagger} G_j^{\rm c}$.
Introducing the complex Wishart matrix $W = G_j^{\rm c \, \dagger} G_j^{\rm c}$, we know that the corresponding Jacobian is proportional to
$(\det A)^{d - k}$ (see e.g.~\cite[Eq.~(3.23)]{Fo10}). Making use too of the simple identity
\begin{equation}\label{uA}
{d \over d\mu}  (\det W)^\mu \Big |_{\mu = 0} = {\rm Tr} \log \det  W
\end{equation}
it is therefore possible to rewrite (\ref{8}) in the case $\Sigma = I_d$ as
\begin{equation}\label{8a}
\mu_1 + \cdots + \mu_k =   {1 \over 2}    {d \over d \mu } \Big \langle  (\det W)^{\mu + d - k}  e^{- {\rm Tr} \, W} \Big \rangle_{W > 0} \Big |_{\mu = 0},
\end{equation} 
where the average is over all positive definite $k \times k$ complex Hermitian matrices.

We see that (\ref{8a}) is a function only of the eigenvalues of $W$. Changing variables to the eigenvalues and eigenvectors (see e.g.~\cite[Proposition 1.3.4]{Fo10}) gives
\begin{equation}\label{8b}
\mu_1 + \cdots + \mu_k =   {1 \over 2 Z_{d, d - k} }    {d \over d \mu } Z_{d,d-k+\mu} \Big |_{\mu = 0},
\end{equation}
where
\begin{equation}\label{Zcd}
Z_{d,c} := \int_0^\infty d x_1 \cdots  \int_0^\infty d x_d \, \prod_{j=1}^d e^{-x_j} x_j^{c} \prod_{1 \le j < l \le d} (x_l - x_j)^2.
\end{equation}
But the integral $Z_{c,k}$ is a particular limiting case of the Selberg integral and as such has a product of gamma function evaluation
(see e.g.~\cite[Prop.~4.7.3 with $\beta = 2$]{Fo10}), telling us that
\begin{align}\label{8c}
\mu_1 + \cdots + \mu_k & =   {1 \over 2 } {d \over d \mu }  \prod_{j=0}^{k-1} {\Gamma( d - k + \mu + 1 + j) \over
\Gamma( d - k + 1 + j)} \nonumber \\
& = {1 \over 2} \sum_{j=0}^{k-1} \Psi( d - k + 1 + j) \: = {1 \over 2} \sum_{j=0}^{k-1} \Psi(d - j) .
\end{align}
The result (\ref{r1}) is now immediate.

Let's now consider Corollary \ref{C1}. We thus want to evaluate (\ref{8.1}) in the case $k=d$ but with $\Sigma$ a general positive definite
matrix. For $d=k$ we can write $ \det   ( G_{k}^{\rm c \, \dagger} \Sigma  G_{k}^{\rm c} )^{1/2}  =  
\det   (   G_{k}^{\rm c} G_{k}^{\rm c \, \dagger} \Sigma)^{1/2} $ since then $ G_{k}^{\rm c} G_{k}^{\rm c \, \dagger} $ has full rank. Furthermore
the probability density corresponding to ${\mathcal N}_{d \times k}^{\rm c}(0,1)$ is then proportional to 
$\exp(-{\rm Tr} \,  G_{k}^{\rm c} G_{k}^{\rm c \, \dagger} )$. Introducing the complex Wishart matrix $W = G_d^{\rm c}G_d^{\rm c \, \dagger} $, for which
the corresponding Jacobian is just a constant, we therefore have
$$
\mu_1 + \cdots + \mu_d = {1 \over 2} \Big ( \log \det \Sigma + 
 \Big \langle  \log (\det W) \,  e^{- {\rm Tr} \, W} \Big \rangle_{W > 0} \Big ).
$$
 Noting that the latter average is just (\ref{8c}) in the case $k=d$ gives (\ref{15}).

\subsection{Proof of Proposition \ref{P2} }
Writing $X = \Sigma^{1/2} G_k^{\rm c}$, (\ref{8.1}) can be rewritten
\begin{equation}\label{8e}
\mu_1 + \cdots + \mu_k =  \Big \langle \log \det  \Big ( X^\dagger X \Big )^{1/2} \Big 
\rangle_{X \in \Sigma^{-1/2} {\mathcal N}_{d \times k}^{\rm c}(0,1)}.
\end{equation}
Now the set of matrices $\Sigma^{-1/2} {\mathcal N}_{d \times k}^{\rm c}(0,1)$ have probability density function proportional to
$\exp ( - {\rm Tr}(X X^\dagger \Sigma^{-1}))$. The $d \times d$ matrix $X X^\dagger$ has rank $k \le d$, and so has $d-k$
zero eigenvalues. This feature, in the case $d < k$, distinguishes the set of matrices  $\Sigma^{-1/2} {\mathcal N}_{d \times k}^{\rm c}(0,1)$
from $ {\mathcal N}_{d \times k}^{\rm c}(0,1) \Sigma_k^{-1/2}$. Thus the latter has probability density proportional to
$\exp ( - {\rm Tr}(X^\dagger X \Sigma_k^{-1}))$ with $X^\dagger X$ of full rank and thus having no zero eigenvalue. These matrices
are termed complex Wishart matrices of mean zero and covariance $\Sigma_k$, while the matrices $X X^\dagger$ are sometimes referred to
as complex pseudo Wishart matrices of mean zero and covariance $\Sigma$.

Let $\rho_{(1)}(\lambda;\Sigma)$ denote the eigenvalue density of the nonzero eigenvalues $\{\lambda_j\}_{j=1,\dots,k}$ for the ensemble
of complex pseudo Wishart matrices so specified. Noting that
\begin{equation}\label{lz}
 \log \det  \Big ( X^\dagger X \Big )^{1/2}  = {1 \over 2} \sum_{j=1}^k \log \lambda_j
 \end{equation}
and so is a linear statistic in $\{\lambda_j\}$, it follows from (\ref{8e})
\begin{equation}\label{lza}
\mu_1 + \cdots + \mu_k =  {1 \over 2} \int_0^\infty ( \log \lambda)  \rho_{(1)}(\lambda;\Sigma) \, d \lambda.
 \end{equation}
 
 Studies in wireless communications \cite{SM04,SMM05,Gh09} have required the same averaged linear statistic, generalized so that
 $\log \lambda \mapsto \log(\lambda - z)$. In these references, methods involving integration over the unitary group have been used to find
 an explicit expression for this generalization of the RHS of (\ref{lza}) in terms of a $d \times d$ matrix. Thus with the eigenvalues of
 $\Sigma^{-1}$ denoted as in Proposition \ref{P2} it is shown
 \begin{equation}\label{L0}
 \int_0^\infty \log (\lambda - z) \rho_{(1)}(\lambda;\Sigma) \, d \lambda =
 C_{k,d}(\{y_i\}) \sum_{m=1}^k \det [ L_m]
 \end{equation}
 where $L_m$ is the $d \times d$ matrix with entries
  \begin{equation}\label{L1}
  (L_m)_{ij} = \left \{
  \begin{array}{ll} \int_0^\infty \log ( t - z) e^{-y_j t} t^{k-i} \, dt, & i = m \\
  y_j^{-(k-i+1)} \Gamma(k-i+1), & i \ne m, \: i \le k \\
  y_j^{d-i}, & i > k
  \end{array} \right.
  \end{equation}
  and 
   \begin{equation}\label{Ckd}
  C_{k,d}(\{y_j\}) = {(-1)^{k((k+1)/2 - d)} \prod_{j=1}^d y_j^k \over
  \prod_{i=1}^{k-1} i! ( \prod_{1 \le i < j \le d} (y_i - y_j) ) }.
  \end{equation}
  
  For purposes of computing (\ref{Lq}), and also for purpose of making the
  presentation more self contained,  it is of interest to revise the derivation of (\ref{L0}). The first step is to consider $d \times k$ matrices $X$ with
  a probability density function 
  $$
  P(X) = {1 \over \pi^{d k /2}} \det \Sigma^{-k } \exp ( - {\rm Tr}(X X^\dagger \Sigma^{-1})).
  $$
  Changing variables  to the nonzero eigenvalues and eigenvectors of $X X^\dagger$ by writing 
  $X X^\dagger = U {\rm diag}\,(\lambda_1,\dots\lambda_k,0,\dots,0) U^\dagger$
  for $U$ unitary shows
  \begin{eqnarray}\label{IZ1}
  P_k(\lambda_1,\dots,\lambda_k) = & \displaystyle
  {   \det \Sigma^{-k}  \over  \prod_{l=0}^{k-1} \Gamma(2+l) \Gamma(d-k+1+l)} \prod_{l=1}^k \lambda_l^{d-k}
  \prod_{1 \le j < l \le k} (\lambda_j - \lambda_l)^2 \nonumber \\
  & \times  \lim_{\lambda_{k+1}, \lambda_{d} \to 0}
\int \exp( - {\rm Tr} (U {\rm diag} (\lambda_1,\dots\lambda_d)U^\dagger  \Sigma^{-1})) \, (U^\dagger dU),
\end{eqnarray}
where $(U^\dagger dU)$ denotes the Haar volume form for $d \times d$ unitary matrices normalized so that
$\int (U^\dagger dU) = 1$. The matrix integral is precisely what is known as the Harish-Chandra--Itzykson-Zuber integral (see
e.g.~\cite[Prop.~11.6.1]{Fo10}), which can be evaluated as a determinant to give 
 \begin{align}\label{IZ1a}
  P_k(\lambda_1,\dots,\lambda_k) = & \displaystyle
  {   \det \Sigma^{-k}  \prod_{l=0}^{d-1} \Gamma(1+l) \over  \prod_{l=0}^{k-1} \Gamma(2+l) \Gamma(d-k+1+l)} \prod_{l=1}^k \lambda_l^{d-k}
  \prod_{1 \le j < l \le  k} (\lambda_j - \lambda_l)^2 \nonumber \\
  &\displaystyle  \times  \lim_{\lambda_{k+1}, \lambda_{d} \to 0}
\prod_{1 \le j < l \le d}{  (y_j - y_l) \over (\lambda_j - \lambda_l)} \det[e^{y_j \lambda_l}]_{j,l=1,\dots,d}.
\end{align} 
The limit can be carried out by power series expanding columns $k+1,\dots,d$, leaving us with the explicit determinant formula
\begin{align}\label{IZ2}
  P_k(\lambda_1,\dots,\lambda_k) = &  \displaystyle
{  (-1)^{k((k+1)/2 - d)} \over \prod_{i=1}^k i!} \prod_{1 \le j < l \le k} (\lambda_j - \lambda_l)
{\prod_{j=1}^d y_j^k \over  \prod_{1 \le j < l \le d} (y_j - y_l)} \nonumber \\
& \times
\det \begin{bmatrix} [ e^{-y_j \lambda_i}]_{i=1,\dots,k \atop j=1,\dots, d} \\
[y_j^{d - i}]_{i=d+1,\dots,k \atop j=1,\dots,d} \end{bmatrix}
\end{align} 

The second step is to use (\ref{IZ2}) to compute the average
 \begin{equation}\label{K0}
\Big  \langle \prod_{j=1}^k (\lambda_j - z)^\mu \Big  \rangle_{P_k}
\end{equation}
For this, one notes that $P_k$ consists of two anti-symmetric factors in $\{\lambda_j\}$. Since, according to the Vandermonde
determinant formula
$$
\prod_{1 \le j < l \le k} (\lambda_j - \lambda_l) =  \mathcal A \, \prod_{l=1}^k \lambda_l^{k-l},
$$
where $ \mathcal A$ denotes the anti-symmetrization operation, we can replace $\prod_{1 \le j < l \le k} (\lambda_j - \lambda_l)$ in the
integrand implied by (\ref{K0}) by $k! \prod_{l=1}^k \lambda_l^{l-1}$. The integrals over $\{\lambda_j\}$ can now be done row-by-row
in the remaining determinant, and we obtain
 \begin{align}\label{K1}
\Big  \langle \prod_{j=1}^k (\lambda_j - z)^\mu \Big \rangle_{P_k} & =
 \displaystyle
{  (-1)^{k((k+1)/2 - d)} \over \prod_{i=1}^{k-1} i!} 
{\prod_{j=1}^d y_j^k \over  \prod_{1 \le j < l \le d} (y_j - y_l)} \nonumber \\
& \times
\det \begin{bmatrix} [ \int_0^\infty (t - z)^\mu t^{k-i} e^{-y_j t} \, dt ]_{i=1,\dots,k \atop j=1,\dots, d} \\
[y_j^{d - i}]_{i=k+1,\dots,d \atop j=1,\dots,d} \end{bmatrix}.
\end{align} 

The third and final step is to differentiate this formula with respect to $\mu$ and set $\mu = 0$. On the LHS this gives the LHS of (\ref{L0}).
Recalling that the differentiation of a determinant with respect to a parameter is equal to the sum of determinants with a single
in each differentiated, we see that the RHS of (\ref{L0}) indeed follows by performing this operation on the RHS of
(\ref{K1}).

Setting $z=0$ in (\ref{L0}) allows the integral in row $i=m$ to be evaluated. Doing this and also taking out appropriate common factors
from each of the first $k$ rows shows
\begin{eqnarray}\label{TA}
\lefteqn{ \int_0^\infty \log  \lambda  \, \rho_{(1)}(\lambda;\Sigma) \, d \lambda} \nonumber \\&&
 = 
{(-1)^{k((k+1)/2 - d)} \over
   \prod_{1 \le i < j \le d} (y_i - y_j)  } 
\sum_{m=1}^k \det \begin{bmatrix} [y_j^{i-1}]_{i=1,\dots,m-1 \atop j =1,\dots,d} \\
[-(\log y_j) y_j^{m-1} +   y_j^{m-1}  \Psi(k-i+1)]_{j=1,\dots,d} \\
[y_j^{d+k-i}]_{i=m+1,\dots,d \atop j =1,\dots,d} \end{bmatrix} .
\end{eqnarray}
Furthermore, reversing the order of the rows $i=k+1, \dots, d$ and recalling that in general a determinant with a single
row having each entry a sum of two terms is equal to the sum of two determinants, it follows from (\ref{TA}) that
\begin{eqnarray}\label{TA}
\lefteqn{ \int_0^\infty \log  \lambda  \, \rho_{(1)}(\lambda;\Sigma) \, d \lambda = 
{1 \over
   \prod_{1 \le i < j \le d} (y_j - y_i)  } \sum_{m=1}^k }  \nonumber \\
&& 
\bigg ( -  \det \begin{bmatrix} [y_j^{i-1}]_{i=1,\dots,m-1 \atop j =1,\dots,d} \\
[(\log y_j) y_j^{m-1} ]_{j=1,\dots,d} \\
[y_j^{d+k-i}]_{i=m+1,\dots,d \atop j =1,\dots,d} \end{bmatrix} +  \Psi(k-m+1) \det[ y_j^{i-1} ]_{i,j=1,\dots,d} \bigg ).   \nonumber \\
\end{eqnarray}
Substituting in (\ref{lza}) and making use of the Vandermonde determinant evaluation gives a result equivalent to (\ref{r2}).

\subsection{Second proof of Corollary \ref{C1}}
Comparison of (\ref{r2}) and (\ref{15}) shows that it suffices to check that
\begin{equation}\label{I1}
 {1 \over  \prod_{1 \le i < j \le d} (y_i - y_j)} \sum_{m=1}^d
\det \begin{bmatrix} [y_j^{i-1}]_{i=1,\dots,m-1 \atop j =1,\dots,d} \\
[(\log y_j) y_j^{m-1}]_{j=1,\dots,d} \\
[y_j^{i-1}]_{i=m+1,\dots,d \atop j =1,\dots,d} \end{bmatrix} =
 \sum_{m=1}^d   \log y_m.
\end{equation}
For this purpose, we observe that the sum over determinants on the LHS is equal to
$$
{d \over d \mu} \det [ y_j^{\mu + i - 1}]_{i,j=1,\dots,d}  \Big |_{\mu = 0}.
$$
But this expression in turn can be rewritten
$$
 \det [ y_j^{ i - 1}]_{i,j=1,\dots,d} {d \over d \mu} \prod_{j=1}^d y_j^\mu \Big |_{\mu = 0} =
  \det [ y_j^{ i - 1}]_{i,j=1,\dots,d }   \sum_{m=1}^d   \log y_m.
 $$
 The identity (\ref{I1}) now follows from the Vandermonde determinant evaluation.
 
 \subsection{The sum rule (\ref{15}) in the case of $A_i = \Sigma^{1/2} G_i^{\rm r}$}
 Suppose the $A_i$ in (\ref{PA}) are given by $A_i = \Sigma^{1/2} G_i^{\rm r}$, where each $G_i^{\rm r}$
 denotes a $d \times d$ real Gaussian with independent standard normals as entries. Let
 $\mathcal N_{d \times d}^{\rm r}(0,1)$ denote the set of real rectangular $d \times k$ Gaussian matrices in
 which each entry is a standard normal. Analogous to (\ref{8.1}) we have
 $$
 \mu_1 + \cdots + \mu_k =   \Big \langle \log \det  \Big ( G_{k}^{\rm r \, T} \Sigma  G_{k}^{\rm r} \Big )^{1/2} \Big 
\rangle_{{\mathcal N}_{d \times k}^{\rm r}(0,1)},
$$
where $G_k^{\rm r}$ denotes $G_k^{\rm r}$ restricted to the first $k$ columns.
As in the argument below (\ref{8c}) in the complex case, for $d=k$ we can write
$ \det   ( G_{k}^{\rm r \, T} \Sigma  G_{k}^{\rm r} )^{1/2}  =  
\det   (   G_{k}^{\rm r} G_{k}^{\rm r \, \dagger} \Sigma)^{1/2} $. Furthermore
the probability density corresponding to ${\mathcal N}_{d \times k}^{\rm r}(0,1)$ is then proportional to 
$\exp(-{\rm Tr} \,  G_{k}^{\rm r} G_{k}^{\rm r \, T}/2 )$. Introducing the real Wishart matrix $W = G_d^{\rm r}G_d^{\rm r \, T} $, for which
the corresponding Jacobian is proportional to $(\det W)^{-1/2}$, we therefore have
$$
\mu_1 + \cdots + \mu_d = {1 \over 2} \Big ( \log \det \Sigma + 
 \Big \langle   (\det W)^{-1/2}  \log (\det W) \,  e^{- {\rm Tr} \, W/2} \Big \rangle_{W > 0} \Big ),
$$
where the average is over positive definite $k \times k$ real symmetric matrices. To evaluate this average we make use of (\ref{uA}) to write
$$
  \Big \langle  \log (\det W) \,  e^{- {\rm Tr} \, W/2} \Big \rangle_{W > 0} =
 {d \over d \mu } \Big \langle  (\det W)^{\mu-1/2}  e^{- {\rm Tr} \, W/2} \Big \rangle_{W > 0} \Big |_{\mu = 0},  
 $$
 Arguing now as in the derivation of (\ref{8b}) shows that the RHS is equal to
 $$
  {1 \over \hat{Z}_{d,-1/2}}  {d \over d \mu }  \hat{Z}_{d,-1/2+ \mu} \Big |_{\mu = 0}
 $$
 where
 $$
  \hat{Z}_{d,c} :=   \int_0^\infty d x_1 \cdots  \int_0^\infty d x_d \, \prod_{j=1}^d e^{-x_j/2} x_j^{c} \prod_{1 \le j < l \le d} |x_l - x_j|.
  $$
  
  Like (\ref{Zcd}), this integral has an evaluation in terms of a product of gamma functions
  (see e.g.~\cite[Prop.~4.7.3 with $\beta = 1$]{Fo10}), giving
  \begin{align*}
 \Big \langle  \log (\det W) \,  e^{- {\rm Tr} \, W/2} \Big \rangle_{W > 0} & =
  {d \over d \mu }2^{d \mu}  \prod_{j=0}^{d-1} {\Gamma( \mu + 1/2 + j/2) \over \Gamma( 1/2 + j/2)}   \Big |_{\mu = 0}, \nonumber \\
  & = d \log 2 + \sum_{j=0}^{d-1} \Psi((j+1)/2).
 \end{align*}
 Consequently, the sought analogue  of (\ref{15}) is
 \begin{equation}\label{2.33}
\mu_1 + \cdots + \mu_d = {1 \over 2} \Big ( \log \det \Sigma +   d \log 2 + \sum_{j=0}^{d-1} \Psi((j+1)/2) \Big ).
\end{equation}

 \subsection{The generalized maximum Lyapunov exponent}
 For the matrix norm in (\ref{Lq}) we take $||P_N|| = \sup_{\vec{x}: |\vec{x}| = 1} |P_N \vec{x}|$. This gives
 $$
 L(q) = {\rm sup}_{\vec{x}: |\vec{x}| = 1} \log \Big \langle  |P_N \vec{x}|^q \Big \rangle.
 $$
 Arguing as in the derivation of (\ref{8.1}), for the $A_i$ in (\ref{PA}) given by (\ref{AG}) we have
   \begin{equation}\label{LK2}
  L(q) =  \log \Big \langle \det   \Big ( G_{1}^{\rm c\, \dagger} \Sigma  G_{1}^{\rm c} \Big )^{q/2} \Big 
\rangle_{{\mathcal N}_{d \times 1}^{\rm c}(0,1)}.
\end{equation}
 Now with $P_k$ the probability density function of the nonzero eigenvalues for matrices $X X^\dagger$, with
 $X$ an element of $\Sigma^{-1/2} {\mathcal N}_{d \times k}^{\rm c}(0,1)$ as in (\ref{IZ1})--(\ref{IZ2}), this can
 be rewritten
 \begin{equation}\label{LK}
 L(q) =  \log \Big \langle \lambda_1^{q/2}  \Big \rangle_{P_1}.
 \end{equation}
 The average in (\ref{LK}) is given by (\ref{K1}) with $k=1$, $z=0$ and $\mu = q/2$, telling us that
  \begin{equation}\label{LK1}
 L(q) = \log \bigg ( \Gamma(1+q/2) {1 \over \prod_{1 \le j < l \le d} (y_j - y_l)}
 \det\begin{bmatrix} [y_j^{-q/2}]_{j=1,\dots,d} \\
 [y_j^{i-1}]_{i=2,\dots,d \atop
 j=1,\dots,d} \end{bmatrix} \bigg ).
 \end{equation} 
 
 In the case that $\Sigma = I_d$, (\ref{LK2}) is a function of $ G_{1}^{\rm c\, \dagger}   G_{1}^{\rm c} $, so a change of
 variables analogous to that used in the derivation of (\ref{8a}) gives
 \begin{equation}\label{LK3} 
 L(q) =  \log \langle y^{q/2 + d - 1} e^{-y} \rangle_{y > 0} = \log {\Gamma(q/2 + d) \over \Gamma(d)}.
 \end{equation}
 
 \section{Lyapunov exponents for diffusing complex matrices}
 The Lyapunov exponents as defined below (\ref{PA1}) relate to the dynamics of the linear system specified by the difference equation
 $\vec{x}_{i+1} = A_{i+1} \vec{x}_i$ for given $\vec{x}_0$. As emphasized in \cite{Ne86}, the continuos counterpart of this setting is
 the matrix stochastic differential equation
 $$
 d X(t) = A X(t) dt + B X(t) d W(t)
 $$
 where $A$ and $B$ are fixed $d \times d$ matrices, and $W(t)$ is a $d \times d$ matrix with complex Brownian entries.
 
 In the case that $A$ and $B$ commute, one has
 $$
 X(t) = \exp \Big ( (A - {1 \over 2} B^2)t + B W(t) \Big ) X(0).
 $$
 Of particular interest is the case $A = {1 \over 2} B^2$, $B = I$ so that $X(t) = \exp(W(t)) X(0)$. Suppose furthermore that
 $W(t) = W_1(t) + i W_2(t)$ with $W_1(t)$ and $W_2(t)$ Hermitian matrices of complex Brownian motions. To specify the
 latter, let GUE${}_d(0,\sigma)$ denote the probability density on $d \times d$ Hermitian matrices $H$ proportional to
 $\exp(- {\rm Tr} \, H^2/\sigma^2)$. We then require that $W_j(1)$ has probability density proportional to
 GUE${}_d(0,\sigma_j)$, for $j=1,2$. An analogous specification of $W(t)$ has been given in
  \cite{Ne86} for $W(t)$ consisting of real Brownian entries, and decomposed as $W(t) = S_1(t) + S_2(t)$, where
  $S_1(t)$ is symmetric, and $S_2(t)$ antisymmetric.
  
Generally the matrix $\exp W(1)$ can be constructed as  
$$
e^{W(1)} = \lim_{m \to \infty} e^{C_m(m)} e^{C_m(m-1)} \cdots e^{C_m(1)}
$$
where $C_m(j) := W(j/m) - W((j-1)/m)$. In the above specification of $W(t)$, $C_m(j)$ is independent of $j$ and distributed as
$C/m^{1/2}$, with $C$ specified as  in (\ref{CA}). Consequently $\exp W(1)$ has the distribution of $A_i$ for $A_i$ as specified by
(\ref{CA}).

With this fact established, the argument leading to (\ref{8.1}) can be used to show that for $\exp W(1)$
$$
\mu_1 + \cdots + \mu_k =
\lim_{m \to \infty} m \Big \langle  \log \det \Big ( E_{d \times k}^T e^{C^\dagger/\sqrt{m}} e^{C/\sqrt{m}} E_{d \times k} \Big )  \Big \rangle,
$$
where the average is over matrices $C = H_1 + i H_2$, with $H_j \in {\rm GUE}_d(0,\sigma_j)$ $(j=1,2)$.
Straightforward expansion in powers of $1/\sqrt{m}$ reduces the RHS to
\begin{eqnarray*}
\lefteqn
{{1 \over 2}  \Big \langle {\rm Tr} \Big (  E_{d \times k}^T (C^\dagger + C)^2 E_{d \times k} - ( E_{d \times k}^T (C^\dagger + C) E_{d \times k})^2 \Big )  \Big \rangle }\nonumber \\
&& =
2\Big \langle   {\rm Tr} \Big (  E_{d \times k}^T H_1^2 E_{d \times k} - ( E_{d \times k}^T H_1 E_{d \times k})^2 \Big )   \Big \rangle\nonumber \\
&& =2 \Big \langle \sum_{i=1}^k \sum_{j=k+1}^d |H_{ij}|^2 \Big \rangle  = \sigma_1^2 \sum_{j=1}^k (d-j).
\end{eqnarray*}
 Consequently we have
 \begin{equation}\label{Cv}
 \mu_k = \sigma_1^2 (d - 2k+1)
\end{equation}
(cf.~\cite[eq.~(15)]{Ne86}). Note that this is independent of $\sigma_2$, and that each Lyapunov exponent vanishes for
$\sigma_1 = 0$, corresponding to $C = i H_2$. This latter point follows from $\exp W(1)$ then being a diffusion on $U(d)$, and so
the modulus of the vectors is unchanged under the corresponding flow.
 
 \section{Discussion}
 \setcounter{equation}{0}
 Consider the case $k=1$ of (\ref{7.1}) and thus  the maximal Lyapunov exponent.
 According to Proposition \ref{P2}, for $A_i$ given by  (\ref{AG}), the exact value of the maximal Lyapunov exponent
 is
 \begin{align}\label{W}
\mu_1 &= -   \displaystyle {1 \over 2 \prod_{1 \le i < j \le d} (y_j - y_i)}
\det \begin{bmatrix} 
[ \log y_j]_{j=1,\dots,d} \\
[y_j^{i-1}]_{i=2,\dots,d \atop j =1,\dots,d} \end{bmatrix} - {1 \over 2} \gamma \nonumber \\
& = - {1 \over 2 } \sum_{j=1}^d  \displaystyle {\log y_j \over \prod_{l=1, l \ne j}^d ( 1 - y_j/y_l)} - {1 \over 2} \gamma, 
\end{align}
 where $\gamma$ denotes Euler's constant. In obtaining the first line the fact that $\Psi(1) = - \gamma$ has been used,
 while the second line follows from the first by expanding the determinant by the first row and using the Vandermonde
 determinant formula.
 Note that replacing $y_i \mapsto \sigma^{-2} y_i$ for each $i=1,\dots,d$ changes $\mu_1$ by
 $\mu_1 \mapsto \mu_1 + \log \sigma$. To see this from the first line in (\ref{W}) requires using the 
 Vandermonde  determinant formula, while in the second line one requires the identity
 $$
  \sum_{j=1}^d  \displaystyle {1 \over \prod_{l=1, l \ne j}^d ( 1 - y_j/y_l)} = 1
  $$
  (see e.g.~\cite[displayed equation below (4.153)]{Fo10}). 
 This change to $\mu_1$ is consistent with the corresponding mapping of the matrices
 $A_i \mapsto \sigma A_i$.
 
 For given distinct $\{y_j\}$ we can use the second equation in (\ref{W}) to give a numerical value of $\gamma_1$.
 For example, with $d=2$, $y_1 = 1$, $y_2 = 1/4$ we obtain $\mu_1 = {4 \over 3} \log 2 - \gamma/2 = 0.63558\dots$. In this
 case
 \begin{equation}\label{AG1}
 A_i = \begin{bmatrix} 1 & 0 \\
 0 & 2 \end{bmatrix} G_i^{\rm c}
 \end{equation}
 where $G_i^{\rm c}$ is a $2 \times 2$ complex Gaussian matrix with entries standard complex normals. With $\vec{x}_0 = [1 \: 0]^T$, let us
 define $\vec{x}_ i = A_i \vec{y}_{i-1}$ $(i=1,2,\dots)$, where $\vec{y}_i := \vec{x}_i/ |\vec{x}_i|$. For a given realization of $\{ A_i \}$, it follows from
 (\ref{7.1}) with $k=1$ that
  \begin{equation}\label{AG2}
 \mu_1 = \lim_{m \to \infty} {1 \over m} \sum_{j=1}^m \log |\vec{x}_j|.
  \end{equation}
 Moreover, straightforward working shows that $ {1 \over m} \sum_{j=1}^m \log |\vec{x}_j|$ has a Gaussian distribution with mean $\mu_1$ and
 a standard deviation $\sigma$ proportional to $1/m^{1/2}$,
 $$
 \sigma^2 = {1 \over m }\bigg ( \Big \langle \Big ( \log \det  \Big ( G_{1}^{\rm c \, \dagger} \Sigma  G_{1}^{\rm c} \Big )^{1/2} \Big )^2 \Big 
\rangle_{{\mathcal N}_{d \times 1}^{\rm c}(0,1)} - \mu_1^2 \bigg )
$$
 This then provides a simple to implement Monte Carlo estimation of $\mu_1$ \cite{CPV93}. In the present setting, with $m=10^6$ we obtained the
 estimation $\mu_1 \approx 0.6341$.
 
 Below (\ref{15}) it was commented that for the average value of the sum of Lyapunov exponents to have a well defined limit for $d \to \infty$ it
 was necessary that the eigenvalues $\{y_m\}$ have the scaling form $y_m/d = Y(y_m/d)$ with $Y(0) = 0$. Under this circumstance it is well
 known (see e.g.~\cite[Th.~7.2.2]{PS11}) that the eigenvalue distribution of $ G^{\rm c \, \dagger} \Sigma  G^{\rm c} $ tends to a well defined nonrandom
 limit with density $u_Y(t)$ say, supported on some interval $I \subset \mathbb R^+$. According to a result of Newman \cite{Ne86a}(see also
 \cite{Ka08}) , one then has
 $\lim_{d \to \infty} e^{\mu_1} = ( \int_I t u_Y(t) dt )^{1/2}$. To derive this from (\ref{W}) does not seem possible, although for a given $Y(x)$
 (\ref{W}) can be used to give a numerical estimation of $\mu_1$. For example, with $Y(x) = 1+x$, computation of (\ref{W}) with
 $d = 5000$ (using high precision arithmetic to avoid catastrophic cancellations) gives $\mu_1 \approx  - 0.183$.
 
 The case $k=d$ of (\ref{7.1}), corresponding to the smallest Lyapunov exponent, admits a form very similar to the second expression
 in (\ref{W}). Thus expanding the determinant by the final row and using the Vandermonde determinant evaluation gives
 \begin{equation}
 \mu_d = - {1 \over 2} \sum_{j=1}^d {\log y_j \over \prod_{l=1, l \ne j} (1 - y_l/y_j)} + {1 \over 2} \Psi(d);
 \end{equation}
 note the interchange of the indices in the denominator relative to (\ref{W}). With $d=2$, $y_1 = 1$, $y_2 = 1/4$ as corresponds to (\ref{AG1}), this
 gives $\mu_2 = - {1 \over 3} \log 2 + {1 \over 2}(1 - \gamma) = - 0.019656\dots$, and thus 
 $\mu_1 + \mu_2 = \log 2 + 1/2  - \gamma = 0.615932\dots$ (note that this is consistent with (\ref{15})).
  Let $\vec{x}_0^{(1)} = [1 \: 0]^T$ and $\vec{x}_0^{(2)} = [0 \: 1]^T$,
 and define $\vec{x}_i^{(p)} = A_i \vec{y}_i^{(p)}$, ($p=1,2$, $i=1,2,\dots$) where $\{\vec{y}_i^{(1)}, \vec{y}_i^{(2)} \}$ is obtained from
 $\{\vec{x}_{i-1}^{(1)}, \vec{x}_{i-1}^{(2)} \}$ by the Gram-Schmidt orthonormalization procedure. For a given realization of $\{A_i\}$, the analogue of
 (\ref{AG2}) is then
 $$
\mu_1 +  \mu_2 = \lim_{m \to \infty}  {1 \over m} \sum_{j=1}^m \log \det (Y_j^\dagger Y_j)^{1/2},
 $$
 where $Y_i$ is the $2 \times 2$ matrix with columns given by $\vec{y}_i^{(1)}$ and  $\vec{y}_i^{(2)}$. This formula without the limit suggests a
 Monte Carlo estimation of $\mu_1 +\mu_2$ \cite{CPV93}. In the present setting, with $A_i$ given by (\ref{AG1}), and choosing $m = 10^6$ gave
 $\mu_1 +  \mu_2 \approx 0.6146$.
 
 \subsection*{Acknowledgements}
 This work was supported by the Australian Research Council.
 
% \bibliographystyle{amsplain}
%\bibliography{book1}

 \providecommand{\bysame}{\leavevmode\hbox to3em{\hrulefill}\thinspace}
\providecommand{\MR}{\relax\ifhmode\unskip\space\fi MR }
\providecommand{\MRhref}[2]{%
  \href{http://www.ams.org/mathscinet-getitem?mr=#1}{#2}
}
\providecommand{\href}[2]{#2}

\end{document}